\documentclass[11pt]{article}
\usepackage{amsmath}
\usepackage{amssymb}
\usepackage{amsbsy}
\usepackage{amsfonts}
\usepackage{mathtools}
\usepackage{bbm}
\usepackage{comment}
\newtheorem{theorem}{Theorem}[section]%
\newtheorem{lemma}[theorem]{Lemma}%
\newtheorem{cor}[theorem]{Corollary}%

\def\f{\noindent}

\newcommand{\qed}{\mbox{\raisebox{0.7ex}{\fbox{}}} \vspace{4truemm}}

\begin{document}

\begin{center}
{\Large{\textbf{On Jacobian group and complexity of the generalized
Petersen graph $GP(n, k)$ through Chebyshev polynomials}}}
\end{center}

\vskip 5mm {\small
\begin{center}
{\textbf{Y.~S.~Kwon,}}\footnote{{\small\em Department of
Mathematics, Yeungnam University, Korea}}
{\textbf{A.~D.~Mednykh,}}\footnote{{\small\em Sobolev Institute of
Mathematics, Novosibirsk State University, Siberian Federal
University}} {\textbf{I.~A.~Mednykh,}}\footnote{{\small\em Sobolev
Institute of Mathematics, Novosibirsk State University, Siberian
Federal University}}
 \end{center}}

\title{ \vspace{-1.2cm}
On Jacobian group and complexity of the generalized Petersen graph
$GP(n, k)$ through Chebyshev polynomials
\thanks{Supported by }}

\begin{abstract}
In the present paper we find a simple algorithm for counting
Jacobian group of the generalized Petersen graph $GP(n, k).$ Also,we
obtain a closed formula for the number of spanning trees of this
graph in terms of Chebyshev polynomials.
\bigskip

\f\textbf{Key Words:} spanning tree,  Jacobian group, Petersen graph, Chebyshev polynomial \\
\textbf{AMS Mathematics Subject Classification:} 05C30, 39A10
\end{abstract}

\section{Introduction}
The notion of the Jacobian group of a graph, which is also known as the Picard group, the critical group, and the dollar
or sandpile group, was independently introduced by many authors (\cite{CoriRoss}, \cite{BakerNorine}, \cite{Biggs},
\cite{BachHarpNagnib}). This notion arises as a discrete version of the Jacobian in the classical theory of Riemann surfaces.
It also admits a natural interpretation in various areas of physics, coding theory, and financial mathematics. The Jacobian group
is an important algebraic invariant of a finite graph. In particular, its order coincides with the number of spanning trees of the
graph, which is known for some simplest graphs, such as the wheel, fan, prism, ladder, and M\"obius ladder [5], grids \cite{NP04},
lattices \cite{SW00}, Sierpinski gaskets \cite{CCY07,AD11}, 3-prism and 3-anti-prism \cite{SWZ16}. At the same time, the structure
of the Jacobian is known only in particular cases \cite{CoriRoss}, \cite{Biggs}, \cite{Lor}, \cite{YaoChinPen}, \cite{ChenHou}, \cite{MedZind}
and \cite{MedMed2}. We mention that the number of spanning trees for circulant graphs is expressed is terms of the Chebyshev polynomials;
it was found in \cite{ZhangYongGol}, \cite{ZhangYongGolin}, and \cite{XiebinLinZhang}. We show that similar results are also true for the
generalized Petersen graph $GP(n, k).$

The generalized Petersen graph $GP(n,k)$ has vertex set and edge set given by
\begin{eqnarray*}
V(P(n,k)) &=& \{ u_i, v_i \ | \  i=1,2, \ldots, n \} \\
E(P(n,k)) &=& \{ u_{i}u_{i+1}, \ u_iv_i, \ v_i v_{i+k}\ | \ i=1,2, \ldots, n \},
\end{eqnarray*}
where the subscripts are expressed as integers modulo $n$. The classical Petersen graph is $P(5,2)$. The number of spanning trees of
the Petersen graph is calculated in \cite{FGW71} and the spectrum of generalized Petersen graphs is obtained in \cite{GS11}. Even
though the number of spanning trees of a regular graph can be computed by its eigenvalues, it is not easy to obtain a closed formula
for the number of spanning trees for $GP(n, k)$ by their result. In this paper we find a closed formula for the number of spanning trees
for $GP(n, k)$ through Chebyshev polynomials. Also, we suggest an effective algorithm for calculating Jacobian of $GP(n, k).$

\section{Basic definitions and preliminary facts}

Consider a connected finite graph $G,$ allowed to have multiple edges but without loops. We endow each edge of $G$
with the two possible directions. Since $G$ has no loops, this operation is well defined. Let $O=O(G)$ be the set
of directed edges of $G.$ Given $e\in O(G),$ we denote its initial and terminal vertices by $s(e)$ and $t(e),$
respectively. Recall that a closed directed path in $G$ is a sequence of directed edges $e_i\in O(G),\, i=1,\ldots,n$
such that $t(e_i)=s(e_{i+1})$  for  $i=1,\,\ldots,n-1$ and $t(e_n)=s(e_1).$

Following \cite{BakerNorine} and \cite{BachHarpNagnib}, \textit{the
Jacobian group}, or simply \textit{Jacobian } $Jac(G)$ of a graph
$G$ is defined as the (maximal) abelian group generated by flows
$\omega(e),e\in O(G),$ obeying the following two Kirchhoff laws:

$K_1:$ the flow through each vertex of $G$ vanishes, that is $\sum\limits_{e\in O,t(e)=x}\omega(e)=0\textrm{ for all }x\in V(G);$

$K_2:$ the flow along each closed directed path $W$ in $G$ vanishes, that is $\sum\limits_{e\in W}\omega(e)=0.$

\noindent Equivalent definitions of the group $Jac(G)$ can be found in papers \cite{CoriRoss}, \cite{BakerNorine}, \cite{Biggs},
\cite{BachHarpNagnib}, \cite{Lor}, \cite{Dhar}, \cite{KotaniSunada}.

We denote the vertex and edge set of $G$ by $V(G)$ and $E(G),$
respectively. Given $u, v\in V(G),$ we set $a_{uv}$ to be equal to
the number of edges between vertices $u$ and $v.$ The matrix $A =
A(G) = \{a_{uv}\}_{u, v\in V (G)},$ called \textit{the adjacency
matrix} of the graph $G.$ The degree $d(v)$ of a vertex $v \in V(G)$
is defined by $d(v) = \sum_u a_{uv}.$ Let $D = D(G)$ be the diagonal
matrix indexed by the elements of $V(G)$ with $d_{vv} = d(v).$
Matrix $L = L(G) = D(G) - A(G)$ is called \textit{the Laplacian
matrix}, or simply \textit{Laplacian}, of the graph $G.$

Recall \cite{Lor} the following useful relation between the
structure of the Laplacian matrix and the Jacobian of a graph $G.$
Consider the Laplacian $L(G)$ as a homomorphism ${\mathbb
Z}^{|V|}\to{\mathbb Z}^{|V|},$ where $|V|=|V(G)|$ is the number of
vertices in $G.$ The cokernel
$\textrm{coker}\,(L(G))=\mathbb{Z}^{|V|}/\textrm{im}\,(L(G))$ --- is
an abelian group. Let
$$\textrm{coker}\,(L(G))\cong\mathbb{Z}_{d_{1}}\oplus\mathbb{Z}_{d_{2}}\oplus\cdots\oplus\mathbb{Z}_{d_{|V|}}$$
be its Smith normal form satisfying the conditions $d_i\big{|}d_{i+1},\,(1\le i\le|V|).$ If the graph is connected,
then the groups ${\mathbb Z}_{d_{1}},{\mathbb Z}_{d_{2}},\ldots,{\mathbb Z}_{d_{|V|-1}}$ --- are finite, and
$\mathbb{Z}_{d_{|V|}}=\mathbb{Z}.$ In this case,
$$Jac(G)\cong\mathbb{Z}_{t_{1}}\oplus\mathbb{Z}_{t_{2}}\oplus\cdots\oplus\mathbb{Z}_{d_{|V|-1}}$$ is the Jacobian
of the graph $G.$ In other words, $Jac(G)$ is isomorphic to the torsion subgroup of the cokernel $\textrm{coker}\,(L(G)).$

Let  $M$ be an integer $n\times n$ matrix, then we can interpret $M$ as a homomorphism from $\mathbb{Z}^n$
to $\mathbb{Z}^n.$ In this interpretation $M$ has a kernel $\textrm{ker}M,$ an image $\textrm{im} M,$ and a
cokernel $\textrm{coker} M = \mathbb{Z}^n/\textrm{im} M.$ We emphasize that $\textrm{coker} M$ of the matrix
$M$ coincides with its Smith normal form.

In what follows, by $I_n$ we denote the identity matrix of order $n.$

We call an $n\times n$ matrix {\it circulant,} and denote it by $circ(a_1, a_2,\ldots,a_n)$ if it is of the form
$$circ(a_0, a_1,\ldots, a_{n-1})=
\left(\begin{array}{ccccc}
a_0 & a_1 & a_2 & \ldots & a_{n-1} \\
a_{n-1} & a_0 & a_1 & \ldots & a_{n-2} \\
  & \vdots &   & \ddots & \vdots \\
a_1 & a_2 & a_3 & \ldots & a_0\\
\end{array}\right).$$

Recall \cite{PJDav} that the eigenvalues of matrix $C=circ(a_0, a_1,\ldots, a_{n-1})$ are given by the following
simple formulas $\lambda_j=p(\varepsilon^j_n),$ where $p(x)=a_0+a_1 x+\ldots+a_{n-1}x^{n-1}$ and $\varepsilon_n$
is the order $n$ primitive root of the unity. Moreover, the circulant matrix $C=p(T),$ where $T=circ(0,1,0,\ldots,0)$
is the matrix representation of the shift operator $T:(x_0, x_1,\ldots,x_{n-2}, x_{n-1})\rightarrow(x_1, x_2,\ldots, x_{n-1},x_0).$
Also, we note that all circulant $n\times n$ matrices share the same set of eigenvectors.

By (\cite{GS11}, lemma 2.1) the $2n\times 2n$ adjacency matrix of the generalized Petersen graph $GP(n,k)$ has the
following block form
$$A(GP(n,k)) =
\left(\begin{array}{cc}
C_n^k & I_n\\
I_n & C_n^1\\
\end{array}\right),$$
where  $C_n^k$ is the $n\times n$ circulant matrix of the form
$C_n^k = circ(\underbrace{0,\ldots,0}_{k\textrm{
times}},1,0,\ldots,0,1,\underbrace{0,\ldots,0}_{k-1\textrm{
times}}).$

Denote by $L=L(GP(n,k))$ the Laplacian of $GP(n,k).$ Since the graph $GP(n,k)$ is three-valent, we have
$$L=3I_{2n}-A(GP(n,k))=\left(\begin{array}{cc}
3I_n-C_n^k & -I_n\\
-I_n & 3I_n- C_n^1\\
\end{array}\right).$$

\smallskip

\section{Cokernels of linear operators given by integer matricies}

Let $P(z)$ be a bimonic integer Laurent polynomial. That is
$P(z)=z^p+a_1z^{p+1}+\ldots+a_{s-1}z^{p+s-1}+z^{p+s}$ for some
integers $p,a_1,a_2,\ldots,a_{s-1}$ and some positive integer $s.$
Introduce the following companion matrix $\mathcal{A}$ for the
polynomial $P(z):$
$\mathcal{A}=\left(\begin{array}{c}\begin{array}{c|c} 0 &
I_{s-1}\end{array}\\\hline -1,-a_1,\ldots,-a_{s-1} \\
\end{array}\right),$ where $I_{s-1}$ is the identity
$(s-1)\times(s-1)$ matrix. We will use the following properties of
$\mathcal{A}.$ Note that $\det \mathcal{A}=(-1)^{s-1}.$ Hence
$\mathcal{A}$ is invertible and inverse matrix $\mathcal{A}^{-1}$ is
also integer matrix. The characteristic polynomial of $\mathcal{A}$
coincides with $z^{-p}P(z).$

Let $R$ be a nonzero commutative ring (integral domain). In most
cases, we deal with the case $R=\mathbb{Z}.$  Denote by
$R^{\mathbb{Z}}$ the set of bi-infinite sequences
$(x_j)_{j\in\mathbb{Z}}=(\ldots,x_{-1},x_0,x_1,x_2,\ldots)$ where
$x_j\in R$ for all $j\in\mathbb{Z}.$ This set is naturally endowed
by the structure of $\mathbb{Z}-$module. Define the shift operator
$T:R^{\mathbb{Z}}\rightarrow R^{\mathbb{Z}}$ by the formula
$T(x_j)_{j\in\mathbb{Z}}=(x_{j+1})_{j\in\mathbb{Z}}.$ For an
arbitrary integer $\ell$ we have
$T^\ell(x_j)_{j\in\mathbb{Z}}=(x_{j+\ell})_{j\in\mathbb{Z}}.$
Simplifying notation we will write $T^\ell x_j=x_{j+\ell}.$ We set
$\mathbbm{1}$ to be the identity operator in $R^\mathbb{Z}.$ For the
sake of simplicity for any integer $n$, we will write $n$ instead of
$n \mathbbm{1}.$ We will use the following notation for the
infinitely generated abelian group. Let $A_\xi,\xi\in \Xi,$ be a
family of $\mathbb{Z}-$linear operators in the space $R^\mathbb{Z}.$
Then by $\langle x \Large| A_\xi x=0, \xi\in \Xi \rangle,$ where
$x=(x_j)_{j\in\mathbb{Z}}$ we denote the abelian group generated by
$x_j,j\in\mathbb{Z}$ satisfying the set of relations $A_\xi
x=0,\xi\in \Xi.$

We will use the following lemma.

\begin{lemma}\label{lemma1} Let $T:R^\mathbb{Z}\to R^\mathbb{Z}$ be the shift operator and $R=\mathbb{Z}.$
Consider two operators $A$ and $B$ given by the formulas $A=P(T), B=Q(T),$ where $P(z)$ and $Q(z)$
are Laurent polynomials with integer coefficients. Then
$$\langle x {\Large|} Ax=0, Bx=0 \rangle \cong\textrm{coker}\,A/\textrm{im}(B|_{\textrm{coker}\,A})
\cong\textrm{coker}_{\textrm{coker}\,A}\,B.$$
\end{lemma}

\textbf{Proof.} Consider $R^{\mathbb{Z}}$ as the abelian group
$\mathbb{Z}^\infty$ of all bi-infnite integer sequences provided
with the natural addition. Then $A$ and $B$ can be considered as
endomorphisms of $\mathbb{Z}^\infty.$ Their images $\textrm{im}\,A$
and $\textrm{im}\,B$ are subgroups in $\mathbb{Z}^\infty.$ Denote by
$\langle\textrm{im}\,A,\textrm{im}\,B\rangle$ the subgroup generated
by elements of $\textrm{im}\,A$ and $\textrm{im}\,B.$ Since $P(z)$
and $Q(z)$ are Laurent polynomials the operators $A=P(T)$ and
$B=Q(T)$ do commute. Hence, subgroup $\textrm{im}\,A$ is invariant
under endomorphism $B.$ Indeed, let $y\in\textrm{im}\,A$ then
$By=B(Ax)=A(Bx)\in\textrm{im}\,A.$ The later means that
$B:\mathbb{Z}^\infty\rightarrow\mathbb{Z}^\infty$ induces an
endomorphism of the group
$\textrm{coker}\,A=\mathbb{Z}^\infty/\textrm{im}\,A.$ We denote this
endomorphism by $B|_{\textrm{coker}\,A}.$ We note that the abelian
group $\langle x {\Large|} Ax=0, Bx=0 \rangle$ is naturally
isomorphic to
$\mathbb{Z}^\infty/\langle\textrm{im}\,A,\textrm{im}\,B\rangle.$ So
we have

\begin{eqnarray*}
&&\mathbb{Z}^\infty/\langle\textrm{im}\,A,\textrm{im}\,B\rangle
\cong(\mathbb{Z}^\infty/\textrm{im}\,A)/\textrm{im}\,(B|_{\textrm{coker}\,A})\\
&&\cong\textrm{coker}\,A/\textrm{im}(B|_{\textrm{coker}\,A})
\cong\textrm{coker}_{\textrm{coker}\,A}\,B.
\end{eqnarray*}

The proof of the lemma is finished.  $\hfill \qed$

\section{Jacobian group for the generalized Petersen graph $GP(n,k)$}
\begin{theorem}\label{theorem0} Let $L=L(GP(n, k))$ be the Laplacian of the generalized Petersen graph
$GP(n, k).$ Then $$\textrm{coker}\,L\cong \textrm{coker}(\mathcal{A}^{n}-I),$$ where $\mathcal{A}$ is
$2(k+1)\times 2(k+1)$ companion matrix for the Laurent polynomial $$(3-z^{-1}-z)(3-z^{-k}-z^k)-1.$$
\end{theorem}

\textbf{Proof.} Let $L$ be the Laplacian matrix of the graph $GP(n,k).$ Then, as it was mentioned above, $L$
is a $2n \times 2n$ matrix of the form
$$L=\left(\begin{array}{cc}3I_n-C_n^k & -I_n\\-I_n & 3I_n-C_n^1\\\end{array}\right),$$
where $C_n^k = circ(\underbrace{0,\ldots,0}_{\textrm{k times}},1,0,\ldots,0,1,\underbrace{0,\ldots,0}_{\textrm{k-1 times}}).$

Consider $L$ as a $\mathbb{Z}-$linear operator
$L:\mathbb{Z}^{2n}\rightarrow\mathbb{Z}^{2n}.$ In this case,
$\textrm{coker}(L)$ is the abelian group generated by elements
$x_1,x_2,\ldots,x_{n},y_1,y_2,\ldots,y_{n}$ satisfying the linear
system of equations
$L(x_1,x_2,\ldots,x_{n},y_1,y_2,\ldots,y_{n})^t=0.$ By the property
mentioned in Section 2, the Jacobian of the graph $GP(n,k)$ is
isomorphic to the finite part of cokernel of the operator $L$. So,
it suffices to show that cokernels of operators
$L:\mathbb{Z}^{2n}\rightarrow\mathbb{Z}^{2n}$ and
$\mathcal{A}^n-I_{2s_k}:\mathbb{Z}^{2s_k}\rightarrow\mathbb{Z}^{2s_k}$
are isomorphic. To study the structure of $\textrm{coker}(L)$ we
consider two bi-infinite sequences of elements
$(x_j)_{j\in\mathbb{Z}}=(\ldots,x_{-1},x_0,x_1,x_2,\ldots)$ and
$(y_j)_{j\in\mathbb{Z}}=(\ldots,y_{-1},y_0,y_1,y_2,\ldots).$ By
circularity of the $n \times n$ blocks of matrix $L$ we have the
following representation for cokernel of $L:$

\begin{eqnarray*}\textrm{coker}(L)&=&\langle x_i,y_i,i\in\mathbb{Z}\large| 3 x_j-x_{j-k}-x_{j+k}-y_j=0,\\
&{}&3 y_j-y_{j-1}-y_{j+1}-x_j=0,  x_{j+n}=x_{j}, y_{j+n}=y_{j}, j\in\mathbb{Z}\rangle.\end{eqnarray*}

Consider the operator $L(T):R^\mathbb{Z}\rightarrow R^\mathbb{Z}$ defined by $L(T)=(3-T-T^{-1})(3-T^k-T^{-k})-1.$
Then by making use of the operator notation we can rewrite the cokernel of $L$ in the following way
\begin{eqnarray*}
\textrm{coker}(L)&=&\langle x,y\Large |(3-T^k-T^{-k})x=y, (3-T-T^{-1})y=x, T^{n}x=x, T^{n}y=y\rangle\\
&=&\langle x\Large |(3-T-T^{-1})(3-T^k-T^{-k})x=x, T^{n}x=x\rangle\\
&=&\langle x\Large |((3-T-T^{-1})(3-T^k-T^{-k})-1)x=0, (T^{n}-1)x=0\rangle \\
&=&\langle x\Large | L(T)x=0, (T^{n}-1)x=0\rangle.\\
\end{eqnarray*}

To finish the proof, we apply  Lemma~\ref{lemma1} to the operators
$A=L(T)$ and $B=Q(T)=T^{n}-1.$

By definition $\textrm{coker}\,A$ is generated by the elements $x_j=e_j+\textrm{im}\,A,j\in\mathbb{Z}.$
Since the Laurent polynomial $P(z)=(3-z-z^{-1})(3-z^k-z^{-k})-1$ is bimonic polynomial it can be represented in the form
$P(z)=z^{-k-1}+a_1z^{-k}+\ldots+a_{2k+1}z^{k}+z^{k+1},$ where $a_1,a_2,\ldots,a_{2k+1}$ are integers. Then the companion
matrix $\mathcal{A}$ is
$\left(\begin{array}{c}
\begin{array}{c|c}0 &I_{2k+1}\end{array}\\ \hline
-1,-a_1,\ldots,-a_{2k+1}
\end{array}\right).$ It is easy to see that $\det \mathcal{A}=(-1)^{2k-1}$ and its inverse $\mathcal{A}^{-1}$ is also integer matrix.

For convenience we set $s=2k+2$ to be  the size of matrix
$\mathcal{A}.$  Let $\mathbb{Z}^s$ be an abelian group generated by
the elements $x_1, x_2,\ldots,x_s$. Note that for any integer $j \in
\mathbb{Z}$,
$(x_{j+1},x_{j+2},\ldots,x_{j+s})^t=\mathcal{A}^{j}(x_{1},x_{2},\ldots,x_{s})^t$,
which implies that  each element $x_j,\,j\in\mathbb{Z}$ can be
uniquely expressed as an integer linear combination of the elements
$x_1, x_2,\ldots,x_s$.

Our present aim is to show that
$\textrm{coker}\,A\cong\mathbb{Z}^s.$  Then we describe the action
of the endomorphism $B|_{\textrm{coker}\,A}$ on the
$\textrm{coker}\,A$. Setting $x=(x_j)_{j\in\mathbb{Z}}$ we can write
$\textrm{coker}\,A=\langle x\Large|Ax=0\rangle.$ So we have the
following representation of $\textrm{coker}\,A$.
\begin{eqnarray*}
&&\textrm{coker}\,A=\langle x\Large|Ax=0\rangle=\\
&&=\langle x_j,j\in\mathbb{Z}\Large| x_{\ell}+a_1x_{\ell+1}+\ldots+a_{s-1}x_{\ell+s-1}+x_{\ell+s}=0,\ell\in\mathbb{Z}\rangle\\
&&=\langle x_j,j\in\mathbb{Z}\Large|(x_{\ell+1},x_{\ell+2},\ldots,x_{\ell+s})^t=
\mathcal{A}(x_{\ell},x_{\ell+1},\ldots,x_{\ell+s-1})^t,\ell\in\mathbb{Z}\rangle\\
&&=\langle x_j,j\in\mathbb{Z}\Large|(x_{\ell+1},x_{\ell+2},\ldots,x_{\ell+s})^t=
\mathcal{A}^{\ell}(x_{1},x_{2},\ldots,x_{s})^t,\ell\in\mathbb{Z}\rangle\\
&&=\langle x_1,x_2,\ldots,x_s{\Large|}\emptyset \rangle\cong\mathbb{Z}^s.
\end{eqnarray*}

Since the operators $A=L(T)$ and $T$ commute, the action
$T|_{\textrm{coker}\,A}:x_j\to x_{j+1},\,j\in\mathbb{Z}$ on the
$\textrm{coker}\,A$ is well defined. Now we describe the action of
$T|_{\textrm{coker}\,A}$ on the set of generators
$x_1,x_2,\ldots,x_s.$ For any $i=1, \ldots, s-1$, we have
$T|_{\textrm{coker}}(x_i)=x_{i+1}$ and $T|_{\textrm{coker}\,A}
(x_s)=x_{s+1}=-x_1-a_1 x_2-\ldots-a_{s-2}x_{s-1}-a_{s-1}x_s$. Hence,
the action of $T|_{\textrm{coker}\,A}$ on the $\textrm{coker}\,A$ is
given by the matrix $\mathcal{A}.$ Considering  $\mathcal{A}$ as an
endomorphism of the $\textrm{coker}\,A,$ we can write
$T|_{\textrm{coker}\,A}=\mathcal{A}.$ Finally,
$B|_{\textrm{coker}\,A}=Q(T|_{\textrm{coker}\,A})=Q(\mathcal{A}).$
Applying Lemma~\ref{lemma1} we finish the proof of the theorem. $\hfill \qed$

\medskip

\begin{cor}\label{corollary1} Sandpile group $\textrm{Jac}(GP(n,k))$ of the generalized Petersen graph $GP(n,k)$ is isomorphic to the torsion subgroup of $\textrm{coker}(\mathcal{A}^{n}-I),$ where
$\mathcal{A}$ is the companion matrix for  the Laurent polynomial $(3-z^{-1}-z)(3-z^{-k}-z^k)-1.$
\end{cor}


The Corollary~\ref{corollary1} gives a simple way to count  Jacobian
group $\textrm{Jac}(GP(n,k))$ for small values of $k$  and
sufficiently large numbers $n.$  The numerical results are presented
in Tables $1, 2$ and $3.$

\section{Counting the number  of spanning trees for the generalized Petersen graph $GP(n,k)$}

\begin{theorem}\label{theorem1}The number of spanning trees of the generalized Petersen graph $GP(n,k)$
is given by the formula $$\tau(GP(n,k))=(-1)^{(n-1)(k-1)}n
\prod_{s=1}^{k}\frac{T_n(w_s)-1}{w_s-1},$$ where $w_s,
s=1,2,\ldots,k$ are roots of the order $k$ algebraic equation
$2T_k(w)-\frac{T_k(w)-1}{w-1}-3=0,$ and $T_k(w)$ is the Chebyshev
polynomial of the first kind.
\end{theorem}

\textbf{Proof.} By the celebrated Kirchhoff theorem, the number of spanning trees $\tau_k(n)$ is equal to the
product of nonzero eigenvalues of the Laplacian of a graph $GP(n,k)$ divided by the number of its vertices $2n.$
To investigate the spectrum of Laplacian matrix we note that matrix $C^k_n=T^{-k}+T^{k},$ where $T=circ(0,1,\ldots,0)$
is the $n \times n$ shift operator. The latter equality easily follows from the identity $T^n=I_n.$  Hence,
$$ L=\left(\begin{array}{cc}
3I_n-T^{-k}-T^{k}  & -I_n\\
-I_n & 3I_n-T^{-1}-T \\
\end{array}\right).$$

The eigenvalues of circulant matrix $T$ are $\varepsilon_n^j,$ where
$\varepsilon_n=e^\frac{2\pi i}{n}.$ Since all eigenvalues of $T$ are
distinct, the matrix $T$ is conjugate to the diagonal matrix
$\mathbb{T}=diag(1,\varepsilon_n,\ldots,\varepsilon_n^{n-1})$, where
diagonal entries of
$diag(1,\varepsilon_n,\ldots,\varepsilon_n^{n-1})$ are
$1,\varepsilon_n,\ldots,\varepsilon_n^{n-1}$. To find spectrum of
$L,$ without loss of generality, one can assume that $T=\mathbb{T}.$
Then the $n \times n$ blocks of $L$ are diagonal matrices. This
essentially simplifies the problem of finding eigenvalues of $L.$
Indeed, let $\lambda$ be an eigenvalue of $L$ and
$(x,y)=(x_1,\ldots,x_n,y_1,\ldots,y_n)$ be the respective
eigenvector. Then we have the following system of equations

$$\left\{\begin{array}{cc}
(3I_n-T^{-k}-T^{k})x-y & =\lambda x\\
-x+(3I_n-T^{-1}-T)y & = \lambda y \\
\end{array}.\right.$$

From here we conclude that $y=(3I_n-T^{-k}-T^{k})x-\lambda x =
(3-\lambda-T^{-k}-T^{k})-1)x$. Substituting $y$ in the second
equation, we have
$((3-\lambda-T^{-1}-T)(3-\lambda-T^{-k}-T^{k})-1)x=0$.

Recall the matrices under consideration are diagonal and the
$(j+1,j+1)$-th entry of $T$ is equal to $\varepsilon_n^{j}.$
Therefore, we have
$((3-\lambda-\varepsilon_n^{-j}-\varepsilon_n^j)(3-\lambda-\varepsilon_n^{-j
k}-\varepsilon_n^{j k})-1)x_{j+1}=0$ and
$y_{j+1}=(3-\lambda-\varepsilon_n^{-j k}-\varepsilon_n^{j
k})x_{j+1}.$

So,  for any $j=0,\ldots, n-1$ the matrix $L$ has two eigenvalues,
say $\lambda_{1,j}$ and $\lambda_{2,j}$ satisfying the quadratic equation $(3-\lambda-\varepsilon_n^{-j}-\varepsilon_n^j)
(3-\lambda-\varepsilon_n^{-jk}-\varepsilon_n^{jk})-1=0.$ The corresponding eigenvectors are $(x,y),$ where $x={\bf e}_{j+1} =(0,\ldots,\underbrace{1}_{(j+1)-th},\ldots, 0)$
and $y=(3-\lambda-T^{-k}-T^{k})\textbf{e}_{j+1}$.
In particular, if $j=0$ for $\lambda_{1,0}, \lambda_{2,0}$ we have $(1-\lambda)(1-\lambda)-1=\lambda(\lambda-2)=0.$
That is, $\lambda_{1,0}=0$ and $\lambda_{2,0}=2.$ Since $\lambda_{1,j}$ and $\lambda_{2,j}$ are roots of the same quadratic
equation, we obtain $\lambda_{1,j}\lambda_{2,j}=P(\varepsilon_n^j),$ where $P(z)=(3-z^{-1}-z)(3-z^{-k}-z^{k})-1.$

Now we have $$\tau_k(n)=\frac{1}{2n}\lambda_{2,0}\prod\limits_{j=1}^{n-1}\lambda_{1,j}\lambda_{2,j}=
\frac{1}{n}\prod\limits_{j=1}^{n-1}\lambda_{1,j}\lambda_{2,j}=\frac{1}{n}\prod\limits_{j=1}^{n-1}P(\varepsilon_n^j).$$
\begin{lemma}\label{lemma2} The following identity holds
$$(3-z^{-1}-z)(3-z^{-k}-z^k)-1=2(w-1)(2T_k(w)-\frac{T_k(w)-1}{w-1}-3),$$
where $T_k(w)$ is the Chebyshev polynomial of the first kind and
$w=\frac{1}{2}(z^{-1}+z).$
\end{lemma}

\textbf{Proof.} Let us substitute $z=e^{i\varphi}.$ It is easy to
see that $w=\frac{1}{2}(z^{-1}+z)=\cos\varphi,$ so we have
$T_k(w)=\cos(k\arccos w)=\cos(k\varphi)$. Then the statement of the
lemma is equivalent to the following elementary identity
$$(3-2\cos\varphi)(3-2\cos(k\varphi))-1=2(\cos\varphi-1)(2\cos(k\varphi)-\frac{\cos(k\varphi)-1}{\cos\varphi-1}-3) \vspace{-1cm}$$
$ \hfill \qed$

By Lemma \ref{lemma2}, $P(z)=2(w-1)h_k(w),$ where
$w=\frac{1}{2}(z+z^{-1})$ and $h_k(w)=2T_k(w)-(T_k(w)-1)/(w-1)-3$ is
the  polynomial of degree $k$. Note that $2(w-1) =
\frac{(z-1)^2}{z}$.  Since $h_k(1)=2T_k(1)-T_k^\prime(1)-3=-1-k^2
\neq 0,$ the Laurent polynomial $P(z)$ has the root $z=1$ with
multiplicity two. Hence, the roots of $P(z)$ are
$1,1,z_1,1/z_1,\ldots,z_k,1/z_k,$ where for all $s=1,\ldots,k,
z_s\neq 1$ and $w_s=\frac{1}{2}(z_s+z_s^{-1})$ is a root of equation
$h_k(w)=0.$ We set $H(z)=\prod\limits_{s=1}^{k}(z-z_s)(z-z_s^{-1}).$
Then $P(z)=\frac{1}{z^{k+1}}(z-1)^2 H(z).$

\medskip

\begin{lemma}\label{lemma5}  Let $H(z)=\prod\limits_{s=1}^{k}(z-z_s)(z-z_s^{-1})$ and  $H(1)\neq 0.$ Then
$$\prod\limits_{j=1}^{n-1}H(\varepsilon_n^j)=\prod\limits_{s=1}^{k}\frac{T_n(w_s)-1}{w_s-1},$$  where
$w_s=\frac12(z_s+z_s^{-1}),\,s=1,\ldots,k$  and $T_n(x)$  is the
Chebyshev polynomial of the first kind.
\end{lemma}

\textbf{Proof.} It is easy to check that $\prod\limits_{j=1}^{n-1}(z-\varepsilon_n^j)=\frac{z^n-1}{z-1}$ if $z\neq1.$
Also we note that $\frac12(z^n+z^{-n})=T_n(\frac12(z+z^{-1})).$  By the substitution $z=e^{i\,\varphi}$  the latter follows from the  evident
identity $\cos(n\varphi)=T_n(\cos\varphi).$  Then we have
\begin{eqnarray*}
\prod\limits_{j=1}^{n-1} H(\varepsilon_n^j) &=& \prod\limits_{j=1}^{n-1}\prod\limits_{s=1}^{k}(\varepsilon_n^j-z_s)(\varepsilon_n^j-z_s^{-1})\\
&=&\prod\limits_{s=1}^{k}\prod\limits_{j=1}^{n-1}(z_s-\varepsilon_n^j)(z_s^{-1}-\varepsilon_n^j)\\
&=&
\prod\limits_{s=1}^{k}\frac{z_s^n-1}{z_s-1}\frac{z_s^{-n}-1}{z_s^{-1}-1}=\prod\limits_{s=1}^{k}\frac{T_n(w_s)-1}{w_s-1}.  \vspace{-1cm}
\end{eqnarray*}  
$\hfill \qed$

Note that
$\prod\limits_{j=1}^{n-1}(1-\varepsilon_n^j)=\lim\limits_{z\to1}\prod\limits_{j=1}^{n-1}
(z-\varepsilon_n^j)= \lim\limits_{z\to1}\frac{z^n-1}{z-1}=n$ and
$\prod\limits_{j=1}^{n-1}\varepsilon_n^{j} = (-1)^{n-1}$.  As a
result, taking into account Lemma~\ref{lemma5}, we obtain
\begin{eqnarray*}
\tau_k(n)&=&\frac{1}{n}\prod\limits_{j=1}^{n-1}P(\varepsilon_n^j)=\frac{1}{n}\prod\limits_{j=1}^{n-1}
\frac{(\varepsilon_n^j-1)^2}{(\varepsilon_n^{j})^{k+1}}H(\varepsilon_n^j)=\frac{(-1)^{(n-1)(k+1)}n^2}{n}
\prod\limits_{j=1}^{n-1}H(\varepsilon_n^j)\\
&=&
(-1)^{(n-1)(k-1)}n\prod\limits_{s=1}^{k}\frac{T_n(w_s)-1}{w_s-1}.   
\end{eqnarray*}  
$\hfill \qed$

\begin{cor}\label{corollary2} $\tau(GP(n,k))=n \left|\prod_{s=1}^{k}U_{n-1}(\sqrt{\frac{1+w_s}{2}})\right|^2,$
where $w_s, s=1,2,\ldots,k$  are the same as in
Theorem~\ref{theorem1} and  $U_{n-1}(w)$ is the Chebyshev polynomial
of the second kind.
\end{cor}

\smallskip
\textbf{Proof.} Follows from the identity
$\frac{T_n(w)-1}{w-1}=U_{n-1}^2(\sqrt{\frac{1+w }{2}}).  \hfill \qed $

\section{Applications and Examples}

\subsection{Prism graph $GP(n,1).$}

As a first consequence of theorem \ref{theorem1} we have the following formula for number of spanning trees of
the $n$-prism graph $GP(n,1):$

$$\tau_1(n)=n(T_n(2)-1).$$

This formula is well known and was independently obtained by many authors. For example, by J. Sedl$\acute{\rm a}$c$\check{\rm e}$k,
J.W. Moon, N. Biggs and others (\cite{BoePro}).

\subsection{Graph $GP(n,2).$}

\begin{theorem}\label{theorem2}The number $\tau_2(n)$ of the spanning trees for the generalized Petersen graph
$GP(n,2)$ is equal to $(-1)^n \frac{1}{20} n (\alpha^2 - 29
\beta^2),$ where the integers $\alpha$ and $\beta$ are given by the
equality: $T_n(\frac{1 + \sqrt{29}}{4}) - 1 =\frac{\alpha +
\beta\sqrt{29}}{4}. $

Moreover, $\tau_2(n) = n a(n)^2,$ where the integer sequence $a(n)$ satisfies the following recursive
relation
\begin{align*}
a(n+4) &= a(n+3) + 3 a(n+2) - a(n+1) - a(n),\\
a(0) &= 0, a(1) = 1, a(2) = 1, a(3)=5.
\end{align*}
Note $a(n)$ is \textit{A192422} sequence in the On - Line Encyclopaedia of Integer Sequences.
\end{theorem}
\textbf{Proof.} Note that $2T_2(w)-\frac{T_2(w)-1}{w-1}-3= 4w^2
-2w-7$, and hence two roots of the equation
$2T_2(w)-\frac{T_2(w)-1}{w-1}-3=0$ are $w_1=\frac{1 + \sqrt{29}}{4}$
and $w_2=\frac{1 - \sqrt{29}}{4}$. By Theorem~\ref{theorem1}, we
obtain

\begin{eqnarray*}\tau_2(n)&=&(-1)^{n-1}n\prod_{s=1}^2\frac{T_n(w_s)-1}{w_s-1}\\&=&\frac{(-1)^{n-1}n(\alpha + \beta \sqrt{29})(\alpha -\beta \sqrt{29})}{16(w_1-1)(w_2-1)}
= \frac{(-1)^nn (\alpha^2 - 29 \beta^2) }{20}.
\end{eqnarray*}

To prove the second statement, by Corollary~\ref{corollary2}, we
have $a(n)=|p(n)|,$ where
$p(n)=U_{n-1}(\theta_1)U_{n-1}(\theta_2),\theta_1=\sqrt{\frac{1+w_1}{2}}$
and $\theta_2=\sqrt{\frac{1+w_2}{2}}.$

Recall that Chebyshev polynomial $U_{n}(\theta)$ satisfies the
recursive relation $U_{n+1}(\theta)-2 \theta
U_n(\theta)+U_{n-1}(\theta)=0$ with initial data $U_0(\theta)=1$ and
$U_1(\theta)=2\theta.$ To find the recursive relation for the
sequence $p(n)$, we will use the following lemma.

\medskip
\begin{lemma}\label{lemma4} Let $P(z)$ and $Q(z)$ are polynomials without multiple roots and let $T u(n)=u(n+1)$ be the shift operator.
Suppose the sequences $u(n)$  and $v(n)$ satisfy the recursive relations $P(T)u(n)=0$  and $Q(T)v(n)=0$ respectively. Then the sequence
$p(n)=u(n)v(n)$ satisfies the recursive relations $R(T)p(n)=0,$ where $R(z)$ is the resultant of polynomials $P(\xi)$ and
$\xi^{\rm deg\,Q} Q(\frac{z}{\xi})$ with respect to $\xi$ and ${\rm deg\,Q}$ is degree of $Q(z).$
\end{lemma}

\textbf{Proof.} Let $\lambda_1,\lambda_2,\ldots,\lambda_s$ and $\mu_1,\mu_2,\ldots,\mu_t$ be distinct roots of polynomials $P(z)$
and $Q(z)$ respectively. Then each solution $u(n)$ of the recursive equation $P(T)u(n)=0$ is a linear combination of the functions
$\lambda_1^n,\lambda_2^n,\ldots,\lambda_s^n,$ while solution $v(n)$ of the recursive equation $Q(T)v(n)=0$ is a linear combination
of $\mu_1^n,\mu_2^n,\ldots,\mu_t^n.$ Hence, their product $p(n)$  is a linear combination of the functions $\lambda_j^n\mu_k^n,\,j=1,2,\ldots,s,\,k=1,2,\ldots,t.$ By definition of resultant, we have $R(\lambda_j \mu_k )=0$ for all $j, k$
and the proof of the lemma follows.  $ \hfill \qed $

\medskip
Now we apply Lemma~\ref{lemma4} to the polynomials
$P(z)=z^2-2\theta_1 z+1$ and $Q(z)=z^2-2\theta_2 z+1$. Now the
resultant of polynomials $P(\xi)$ and $\xi^2 Q(\frac{z}{\xi})$ is
\begin{eqnarray*} & \left(\theta_1 +\sqrt{\theta_1^2 -1} -(\theta_2 +\sqrt{\theta_2^2
-1})z\right)\left(\theta_1 +\sqrt{\theta_1^2 -1} -(\theta_2
-\sqrt{\theta_2^2 -1})z\right) \times \\ & \left(\theta_1
-\sqrt{\theta_1^2 -1} -(\theta_2 +\sqrt{\theta_2^2
-1})z\right)\left(\theta_1 -\sqrt{\theta_1^2 -1} -(\theta_2
-\sqrt{\theta_2^2 -1})z\right).
\end{eqnarray*}

If we expand the equation, we get the equation $z^4 - iz^3 +3 z^2 -i
z +1$. Hence we have the following recursive relation for
$p(n)=U_{n-1}(\theta_1)U_{n-1}(\theta_2):$
$$ p(n)-i\, p(n+1)+3p(n+2)-i\, p(n+3)+p(n+4)=0.$$

 Observing that $a(n)=|p(n)|=(-i)^{n-1}p(n)$  and the initial values
$a(0) = 0, a(1) = 1, a(2) = 1, a(3)=5$,   we get the result. $ \hfill \qed $

\subsection{Graph $GP(n,3).$}

\begin{theorem}\label{theorem3}The number $\tau_3(n)$ of the spanning trees for the generalized Petersen graph
$GP(n, 3)$ is given by the formula
$$\tau_3(n)=n\prod_{s=1}^{3}\frac{T_n(w_s)-1}{w_s-1},$$
and $w_1, w_2, w_3$  are roots of the equation  $4w^3-2w^2-5w-2=0.$

Moreover, $\tau_3(2n) = 12 n \,a(n)^2 $ and $\tau_3(2n+1) = (2n+1) b(n)^2,$ where the integer sequences
$a(n)$ and $b(n)$ satisfy the recursive relation
\begin{align*} u(n)&-4u(n+1)-u(n+2)-24u(n+3)+65u(n+4)\\&-24u(n+5)-u(n+6)-4u(n+7)+u(n+8)=0
\end{align*} with the following initial data
\begin{align*} a(0) &=0,a(1)=1, a(2)=4, a(3)=9, a(4)=72,\\ a(5)&=320,u(6)=1332,a(7)=6889
\end{align*} and
\begin{align*} b(0) &=1,b(1)=1, b(2)=20, b(3)=83, b(4) =289,  \\  b(5)&=1693, b(6) =7775,b(7) =34820.
\end{align*}
\end{theorem}

\textbf{Proof.} The first statement of the theorem directly follows
from Theorem~\ref{theorem1}. To prove the second, by
Corollary~\ref{corollary2}, we have $\tau_3(n)=n|c(n)|^2,$ where
$c(n)=\prod_{s=1}^3U_{n-1}(\theta_s)$ and
$\theta_s=\sqrt{\frac{1+w_s}{2}},s=1,2,3.$ Then the sequence
$u(n)=U_{n-1}(\theta_s)$ satisfies the following recursive relation
$P_s(T)u(n)=0,$ where $P_s(z)=z^2-2\theta_s z+1.$ Suppose that
$\lambda,\mu,\nu$ are roots of the equations
$P_1(\lambda)=0,P_2(\mu)=0,P_3(\nu)=0.$ By applying Lemma
\ref{lemma4} twice, we obtain that $\eta=\lambda\,\mu\,\nu$ is a
root of the equation $$1 + \eta^2 + 11 \eta^4 + \eta^6 + \eta^8 =
\sqrt{6}\eta(1 + 2\eta^2 + 2 \eta^4 + \eta^6).$$ Denote by
$\eta_1,\eta_2,\ldots,\eta_8$ distinct roots of the previous
equation. Then the sequence
$c(n)=U_{n-1}(\theta_1)U_{n-1}(\theta_2)U_{n-1}(\theta_3)$ can be
written in the form $c(n)=\sum\limits_{j=1}^{8}c_j\eta_j^n,$ where
$c_1,c_2,\ldots,c_8$ are  suitable constants. Let
$\zeta_j=\eta_j^2$. Now $c(2n)=\sum\limits_{j=1}^{8}c_j\zeta_j^{n}$
and $c(2n+1)=\sum\limits_{j=1}^{8}(c_j \eta_j)\zeta_j^{n}.$ Since
$$(1 + \eta^2 + 11 \eta^4 + \eta^6 + \eta^8)^2 = 6\eta^2(1 + 2\eta^2 + 2 \eta^4 + \eta^6)^2,$$ all $\zeta_j,j=1,2,\ldots,8$ are
the roots of the equation
$$(1 + \zeta + 11 \zeta^2+ \zeta^3+ \zeta^4)^2 - 6\zeta(1 + 2\zeta + 2 \zeta^2 + \zeta^3)^2=0.$$
Hence,
$$1 - 4 \zeta - \zeta^2 - 24 \zeta^3 + 65 \zeta^4 - 24 \zeta^5 - \zeta^6 - 4 \zeta^7 + \zeta^8 = 0$$
and both sequences $c(2n)$ and $c(2n+1)$ are solution of the difference equation
\begin{align*}u(n)&-4u(n+1)-u(n+2)-24u(n+3)+65u(n+4)\\
&-24u(n + 5)-u(n + 6)-4u(n + 7) + u(n + 8)=0.\end{align*}
We use the formulas $a(n)=c(2n)/\sqrt{6}$ and
$b(n)= c(2n+1)$ to  calculate the initial
elements of   sequences $a(n)$ and   $b(n)$ directly.
$ \hfill \qed $

\subsection{Graph $GP(n,4).$}

\begin{theorem}\label{theorem4}The number $\tau_4(n)$ of the spanning trees for the generalized Petersen graph
$GP(n,4)$ is given by the formula
$$\tau_4(n)=(-1)^{n-1}n\prod_{s=1}^{4}\frac{T_n(w_s)-1}{w_s-1},$$ and
$w_1, w_2, w_3, w_4$ are roots of the equation
$16w^4-8w^3-24w^2-1=0.$

Moreover, $\tau_4(n)= n\,a(n)^2,$ where the integer sequences $a(n)$ satisfies the recursive relation
$P(T)a(n)=0,$ where

\begin{align*}
P(T) &= T^{16}-T^{15}-2 T^{13}-16 T^{12}+10 T^{11}-2 T^{10}+16 T^{9} \\
&+50 T^{8}-16 T^{7}-2 T^{6}-10 T^{5}-16 T^{4}+2 T^{3}+T+1
\end{align*}
and $T a(n)=a(n+1)$ is the shift operator.

The initial data of $a(n)$ for $n$ equal to $$-7,-6,-5,-4,-3,-2,-1,0,1,2,3,4,5,6,7,8$$ are, respectively,
$$-83,35,-19,1,-5,1,-1,0,1,1,5,1,19,35,83,73.$$
\end{theorem}
\textbf{Proof.} The first statement of the theorem  directly follows from Theorem~\ref{theorem1}. Applying Lemma~\ref{lemma4} several times and using the same arguments as in the proof of Theorem~\ref{theorem2}
we conclude that the sequence $a(n)=|p(n)|=(-i)^{n-1}p(n),$  where $p(n)=\prod_{s=1}^4U_{n-1}(\sqrt{\frac{1+x_s}{4}})$ satisfies the recursive  relation $P(T)a(n)=0.$ $ \hfill \qed $

\section{Final Remarks and Tables}

Theorem~\ref{theorem0}  is the first step to understand the
structure of the Jacobian for $GP(n,k).$ Also, it gives a simple way
for numerical calculations of $\textrm{Jac}(GP(n,k))$ for small
values of $n$ and $k.$ See Tables $1,2,3$ below. Theorems
\ref{theorem2}, \ref{theorem3} and \ref{theorem4} contain very
convenient formulas for counting the number of spanning trees. The
corresponding numerical results are given in Tables $1, 2$ and $3.$

\begin{table}[h]\label{table1}
\caption{Graph $GP(n, 2)$}
\begin{tabular}{r|l|r}
$n$ & $\textrm{Jac}(GP(n,2))$ &$\tau_2(n) =|\textrm{Jac}(GP(n,2))|$\\ \hline
3  & $\mathbb{Z}_{5}\oplus\mathbb{Z}_{15}$ & 75\\
4  & $\mathbb{Z}_{7}\oplus\mathbb{Z}_{28}$ & 196\\
5  & $\mathbb{Z}_{2}\oplus\mathbb{Z}_{10}\oplus\mathbb{Z}_{10}\oplus\mathbb{Z}_{10}$ & 2000\\
6  & $\mathbb{Z}_{35}\oplus\mathbb{Z}_{210}$ & 7350\\
7  & $\mathbb{Z}_{83}\oplus\mathbb{Z}_{581}$ & 48223\\
8  & $\mathbb{Z}_{161}\oplus\mathbb{Z}_{1288}$ & 207368\\
9  & $\mathbb{Z}_{355}\oplus\mathbb{Z}_{3195}$ & 1134225\\
10 & $\mathbb{Z}_{2}\oplus\mathbb{Z}_{12}\oplus\mathbb{Z}_{60}\oplus\mathbb{Z}_{60}\oplus\mathbb{Z}_{60}$ & 5184000\\
11 & $\mathbb{Z}_{1541}\oplus\mathbb{Z}_{16951}$ & 26121491\\
12 & $\mathbb{Z}_{7}\oplus\mathbb{Z}_{7}\oplus\mathbb{Z}_{1365} \oplus\mathbb{Z}_{1820} $ & 121730700\\
13 & $\mathbb{Z}_{6733}\oplus\mathbb{Z}_{87529}$ & 583332757\\
14 & $\mathbb{Z}_{14027}\oplus\mathbb{Z}_{196378}$ & 2754594206\\
15 & $\mathbb{Z}_{5}\oplus\mathbb{Z}_{10}\oplus\mathbb{Z}_{10}\oplus\mathbb{Z}_{2950}\oplus\mathbb{Z}_{8850}$ & 13053750000\\
16 & $\mathbb{Z}_{61663}\oplus\mathbb{Z}_{986608}$ & 60837209104\\
17 & $\mathbb{Z}_{129403}\oplus\mathbb{Z}_{2199851}$ & 284667318953\\
18 & $\mathbb{Z}_{270865}\oplus\mathbb{Z}_{4875570}$ & 1320621268050\\
19 & $\mathbb{Z}_{567911}\oplus\mathbb{Z}_{10790309}$ & 6127935174499\\
20 & $\mathbb{Z}_{4}\oplus\mathbb{Z}_{24}\oplus\mathbb{Z}_{120}\oplus\mathbb{Z}_{49560}\oplus\mathbb{Z}_{49560}$ & 28295350272000\\
\end{tabular}
\end{table}

\clearpage

\begin{table}[h]\caption{Graph $GP(n, 3)$} 
 \begin{tabular}{r|l|r}  $n$ & $\textrm{Jac}(GP(n,3))$
&$\tau_3(n)=|\textrm{Jac}(GP(n,3))|$ \\ \hline
4  & $\mathbb{Z}_{2}\oplus\mathbb{Z}_{8}\oplus\mathbb{Z}_{24}$ & 384\\
5  & $\mathbb{Z}_{2}\oplus\mathbb{Z}_{10}\oplus\mathbb{Z}_{10}\oplus\mathbb{Z}_{10}$ & 2000\\
6  & $\mathbb{Z}_{3}\oplus\mathbb{Z}_{9}\oplus\mathbb{Z}_{108}$ & 2916\\
7  & $\mathbb{Z}_{83}\oplus\mathbb{Z}_{581}$ & 48223\\
8  & $\mathbb{Z}_{3}\oplus\mathbb{Z}_{3}\oplus\mathbb{Z}_{12}\oplus\mathbb{Z}_{48}\oplus\mathbb{Z}_{48}$ & 248832\\
9  & $\mathbb{Z}_{289}\oplus\mathbb{Z}_{2601}$ & 751689\\
10 & $\mathbb{Z}_{4}\oplus\mathbb{Z}_{8}\oplus\mathbb{Z}_{40}\oplus\mathbb{Z}_{40}\oplus\mathbb{Z}_{120}$ & 6144000\\
11 & $\mathbb{Z}_{1693}\oplus\mathbb{Z}_{18623}$ & 31528739\\
12 & $\mathbb{Z}_{6}\oplus\mathbb{Z}_{2664}\oplus\mathbb{Z}_{7992}$ & 127744128\\
13 & $\mathbb{Z}_{5}\oplus\mathbb{Z}_{5}\oplus\mathbb{Z}_{1555}\oplus\mathbb{Z}_{20215}$ & 785858125\\
14 & $\mathbb{Z}_{83}\oplus\mathbb{Z}_{83}\oplus\mathbb{Z}_{83}\oplus\mathbb{Z}_{6972}$ & 3986498964\\
15 & $\mathbb{Z}_{2}\oplus\mathbb{Z}_{10}\oplus\mathbb{Z}_{52230}\oplus\mathbb{Z}_{17410}$ & 18186486000\\
16 & $\mathbb{Z}_{3}\oplus\mathbb{Z}_{3}\oplus\mathbb{Z}_{24}\oplus\mathbb{Z}_{21408}\oplus\mathbb{Z}_{21408}$ & 98993332224\\
17 & $\mathbb{Z}_{170917}\oplus\mathbb{Z}_{2905589}$ & 496614555113\\
18 & $\mathbb{Z}_{9}\oplus\mathbb{Z}_{148257}\oplus\mathbb{Z}_{1779084}$ & 2373854909292\\
19 & $\mathbb{Z}_{802141}\oplus\mathbb{Z}_{15240679}$ & 12225173493739\\
20 & $\mathbb{Z}_{2}\oplus\mathbb{Z}_{8}\oplus\mathbb{Z}_{8}\oplus\mathbb{Z}_{16}\oplus\mathbb{Z}_{80}\oplus\mathbb{Z}_{11120}\oplus\mathbb{Z}_{33360}$ & 60778610688000  \\
\end{tabular} \vspace{1cm}
\caption{ Graph $GP(n, 4)$}  
\begin{tabular}{r|l|r}
$n$ & $\textrm{Jac}(GP(n,4))$ &$\tau_4(n)=|\textrm{Jac}(GP(n,4))|$ \\ \hline
5 & $\mathbb{Z}_{19}\oplus\mathbb{Z}_{95}$ & 1805\\
6 & $\mathbb{Z}_{35}\oplus\mathbb{Z}_{210}$ & 7350\\
7 & $\mathbb{Z}_{83}\oplus\mathbb{Z}_{581}$ & 48223\\
8 & $\mathbb{Z}_{73}\oplus\mathbb{Z}_{584}$ & 42632\\
9 & $\mathbb{Z}_{355}\oplus\mathbb{Z}_{3195}$ & 1134225\\
10& $\mathbb{Z}_{779}\oplus\mathbb{Z}_{7790}$ & 6068410\\
11& $\mathbb{Z}_{1693}\oplus\mathbb{Z}_{18623}$ & 31528739\\
12& $\mathbb{Z}_{2555}\oplus\mathbb{Z}_{30660}$ & 78336300\\
13& $\mathbb{Z}_{5}\oplus\mathbb{Z}_{5}\oplus\mathbb{Z}_{1555}\oplus\mathbb{Z}_{20215}$ & 785858125\\
14& $\mathbb{Z}_{17513}\oplus\mathbb{Z}_{245182}$ & 4293872366\\
15& $\mathbb{Z}_{2}\oplus\mathbb{Z}_{2}\oplus\mathbb{Z}_{2}\oplus\mathbb{Z}_{10}\oplus\mathbb{Z}_{10}
\oplus\mathbb{Z}_{10}\oplus\mathbb{Z}_{950}\oplus\mathbb{Z}_{2850}$ & 21660000000\\
16& $\mathbb{Z}_{71321}\oplus\mathbb{Z}_{1141136}$ & 81386960656\\
17& $\mathbb{Z}_{103}\oplus\mathbb{Z}_{1751}\oplus\mathbb{Z}_{1751}\oplus\mathbb{Z}_{1751}$ & 552962478353\\
18& $\mathbb{Z}_{405055}\oplus\mathbb{Z}_{7290990}$ & 2953251954450\\
19& $\mathbb{Z}_{37}\oplus\mathbb{Z}_{37}\oplus\mathbb{Z}_{23939}\oplus\mathbb{Z}_{454841}$ & 14906272578931\\
20& $\mathbb{Z}_{1823639}\oplus\mathbb{Z}_{36472780}$ & 66513184046420\\
\end{tabular} 
\end{table}

\clearpage

\section*{ACKNOWLEDGMENTS}
The first author was supported by the 2015 Yeungnam University
Research Grant. The second and the third authors were  supported by
the Grant of the Russian Federation Govement at Siberian Federal
University (grant no. 14.Y26.31.0006), by the Presidium of the
Russian Academy of Sciences (project no. 0314-2015-0011), and by the
Russian Foundation for Basic Research (projects no. 15-01-07906 and
16-31-00138).

\end{document}